\documentstyle[fleqn,12pt]{article}
\begin{document}\pagenumbering{arabic}\setcounter{page}{1}\pagestyle{plain}\baselineskip=18pt

\thispagestyle{empty} \rightline{YTUMB 2006-02, July 2006}
\vspace{1.4cm}

\begin{center}
{\Large\bf CARTAN CALCULUS ON THE QUANTUM SPACE ${\cal R}_q^{3}$ }
\end{center}

\vspace{1cm}\noindent Salih \c Celik$^{1,2}$, E. Mehmet \"Ozkan$^1$ and Erg\"un Ya\c sar$^1$

\vspace{0.3cm}\noindent $^1$ Yildiz Technical University, Department of Mathematics, 34210 DAVUTPASA-Esenler, Istanbul, TURKEY.

\vspace{0.3cm}\noindent $^2$ E-mail: sacelik@yildiz.edu.tr

\vspace{3cm} {\bf ABSTRACT}

\vspace{0.3cm}\noindent To give a Cartan calculus on the extended quantum 3d space, the noncommutative differential calculus on the extended quantum 3d space is extended by introducing inner derivations and Lie derivatives.

\newpage\noindent
{\bf 1. INTRODUCTION}

\vspace*{0.3cm}\noindent The noncommutative differential geometry of quantum groups was introduced by Woronowicz [11,12]. In this approach the differential calculus on the group is deduced from the properties of the group and it involves functions on the group, differentials, differential forms and derivatives. The other approach, initiated by Wess and Zumino [10], followed Manin's emphasis [5] on the quantum spaces as the primary objects. Differential forms are defined in terms of noncommuting coordinates, and the differential and algebraic properties of quantum groups acting on these spaces are obtained from the properties of the spaces.

\vspace*{0.3cm}\noindent The differential calculus on the quantum 3d space similarly involves functions on the 3d space, differentials, differential forms and derivatives. The exterior derivative is a linear operator {\sf d} acting on $k$-forms and producing $(k+1)$-forms, such that for scalar functions (0-forms) $f$ and $g$ we have
\begin{eqnarray*}
 {\sf d} (1) & = & 0, \\
 {\sf d}(f g) & = & ({\sf d} f) g + (-1)^{deg(f)} \, f \, ({\sf d}g)
\end{eqnarray*}
where $deg(f) = 0$ for even variables and $deg(f) = 1$ for odd variables, and for a $k$-form $\omega_1$ and any form $\omega_2$
\begin{eqnarray*}
    {\sf d}(\omega_1 \wedge \omega_2) & = & ({\sf d} \omega_1) \wedge \omega_2+(-1)^k \, \omega_1 \, \wedge ({\sf d}\omega_2).
\end{eqnarray*}
A fundamental property of the exterior derivative {\sf d} is
\begin{eqnarray*}
    {\sf d} \wedge {\sf d} & = & : {\sf d}^2  =  0.
\end{eqnarray*}
\noindent There is a relationship of the exterior derivative with the Lie derivative and to describe this relation, we introduce a new operator: the inner derivation. Hence the differential calculus on the quantum 3d space can be extended into a large calculus. We call this new calculus the Cartan calculus. The connection of the inner derivation denoted by {\bf \textit i}$_a$ and the Lie derivative denoted by ${\cal L}_a$ is given by the Cartan formula:
\begin{eqnarray*}
{\cal L}_a & = & {\bf \textit i}_a \circ {\sf d} + {\sf d} \circ {\bf \textit i}_a.
\end{eqnarray*}
This and other formulae are explaned in Ref. 6-8. We now shall give a brief overview without much discussion.

\vspace*{0.3cm} \noindent Let us begin with some information about the inner derivations. Generally, for a smooth vector field $X$ on a manifold the inner derivation, denoted by ${\bf \textit i}_X$, is a linear operator which maps $k$-forms to $(k-1)$-forms. If we define the inner derivation ${\bf \textit i}_X$ on the set of all differential forms on a manifold, we know that ${\bf \textit i}_X$ is an antiderivation of degree $- 1$:
\begin{eqnarray*}
    {\bf \textit i}_X (\alpha \wedge \beta) & = & ({\bf \textit i}_X \alpha) \wedge \beta + (-1)^k \, \alpha \wedge ({\bf \textit i}_X \beta)
\end{eqnarray*}
where $\alpha$ and $\beta$ are both differential forms. The inner derivation ${\bf \textit i}_X$ acts on 0- and 1-forms as follows:
\begin{eqnarray*}
    {\bf \textit i}_X (f) & = & 0, \\
    {\bf \textit i}_X ({\sf d} f) & = & X(f).
\end{eqnarray*}

\vspace*{0.3cm} \noindent We know, from the classical differential geometry, that the Lie derivative ${\cal L}$ can be defined as a linear map from the exterior algebra into itself which takes $k$-forms to $k$-forms. For a 0-form, that is, an ordinary function $f$, the Lie derivative is just the contraction of the exterior derivative with the vector field $X$:
\begin{eqnarray*}
    {\cal L}_X f & = & {\bf \textit i}_X \, {\sf d} f.
\end{eqnarray*}
For a general differential form, the Lie derivative is likewise a contraction, taking into account the variation in $X$:
\begin{eqnarray*}
    {\cal L}_X \, \alpha & = & {\bf \textit i}_X \, {\sf d} \alpha + {\sf d} ({\bf \textit i}_X \alpha).
\end{eqnarray*}
The Lie derivative has the following properties. If ${\cal F}(M)$ is the algebra of functions defined on the manifold $M$ then
\begin{eqnarray*}
{\cal L}_X : {\cal F}(M) \longrightarrow {\cal F}(M)
\end{eqnarray*}
is a derivation on the algebra ${\cal F}(M)$:
\begin{eqnarray*}
    {\cal L}_X (a f + b g) & = & a ({\cal L}_X f) + b ({\cal L}_X g), \\
  {\cal L}_X (f g) & = & ({\cal L}_X f) \, g + f \, ({\cal L}_X g),
\end{eqnarray*}
where $a$ and $b$ real numbers.

\vspace*{0.3cm} \noindent The Lie derivative is a derivation on ${\cal F}(M) \times {\cal V}(M)$ where ${\cal V}(M)$ is the set of vector fields on $M$:
\begin{eqnarray*}
    {\cal L}_{X_1} (f X_2) & = & ({\cal L}_{X_1} f) \, X_2 + f \, ({\cal L}_{X_1} X_2).
\end{eqnarray*}
The Lie derivative also has an important property when acting on differential forms. If $\alpha$ and $\beta$ are two differential forms on $M$ then
\begin{eqnarray*}
    {\cal L}_X (\alpha \wedge \beta) & = & ({\cal L}_X \alpha) \wedge \beta + (-1)^k \, \alpha \wedge ({\cal L}_X \beta)
\end{eqnarray*}
where $\alpha$ is a $k$-form.

\vspace*{0.3cm}\noindent The extended calculus on the quantum plane was introduced in Ref. 3 using the approach of Ref. 6. In this work we explicitly set up the Cartan calculus on the quantum 3d space using approach of Ref 1.

\vspace*{0.5cm} \noindent {\bf 2. REVIEW OF SOME STRUCTURES ON ${\cal R}^3_q$}

\vspace*{0.3cm}\noindent In this section we give some information on the Hopf algebra structures of the quantum 3d space and its differential calculus [2] which we shall use in order to establish our notions.

\vspace*{0.3cm} \noindent {\bf 2.1 The algebra of polynomials on the quantum 3d space}

\vspace*{0.3cm} \noindent The quantum three dimensional space is defined as an associative algebra generated by three noncommuting coordinates $x$, $y$ and $z$ with three quadratic relations
\begin{eqnarray*}
x y & = & q y x, \\ y z & = & q z y, \\ x z & = & q z x,
\end{eqnarray*}
where $q$ is a non-zero complex number. This associative algebra over the complex number, C, is known as the algebra of polynomials over the quantum three dimensional space and we shall denote it by ${\cal R}_q^{3 }$. In the limit $q \longrightarrow 1$, this algebra is commutative and can be considered as the algebra of polynomials C$[x,y,z]$ over the usual three dimensional space, where $x$, $y$ and $z$ are the three coordinate functions. We denote the unital extension of ${\cal R}_q^{3}$ by ${\cal A}$, i.e. it is obtained by adding a unit element.

\vspace*{0.3cm} \noindent {\bf 2.2 The Hopf algebra structure on ${\cal A}$}

\vspace*{0.3cm} \noindent One extends the algebra ${\cal A}$ by including inverse of $x$ which obeys
\begin{eqnarray*}
    x x^{-1} = 1 =  x^{-1} x.
\end{eqnarray*}

\vspace*{0.3cm} \noindent The definitions of a coproduct, a counit and a coinverse on the algebra ${\cal A}$ as follows [2]:

\vspace*{0.3cm} {\bf (1)} The C-algebra homomorphism (coproduct) $\Delta_{\cal A}: {\cal A} \longrightarrow {\cal A} \otimes {\cal A}$ is defined by
\begin{eqnarray*}
\Delta_{\cal A}(x) & = & x \otimes x, \\
\Delta_{\cal A}(y) & = & x \otimes y + y \otimes x,  \\
\Delta_{\cal A}(z) & = & z \otimes 1 + 1 \otimes z,
\end{eqnarray*}
which is coassociative:
\begin{eqnarray*}
   (\Delta_{\cal A} \otimes \mbox{id}) \circ \Delta_{\cal A} = (\mbox{id} \otimes \Delta_{\cal A}) \circ \Delta_{\cal A}
\end{eqnarray*}
where id denotes the identity map on ${\cal A}$.

\vspace*{0.3cm} {\bf (2)} The C-algebra homomorphism (counit) $\epsilon_{\cal A}: {\cal A} \longrightarrow {\mbox C}$ is given by
\begin{eqnarray*}
    \epsilon_{\cal A}(x) & = & 1, \\ \epsilon_{\cal A}(y) & = & 0, \\ \epsilon_{\cal A}(z) & = & 0.
\end{eqnarray*}
The counit $\epsilon_{\cal A}$ has the property
\begin{eqnarray*}
    \mu \circ (\epsilon_{\cal A} \otimes \mbox{id}) \circ \Delta_{\cal A} = \mu' \circ (\mbox{id} \otimes \epsilon_{\cal A}) \circ \Delta_{\cal A}
\end{eqnarray*}
where $\mu:{\mbox C}\otimes{\cal A}\longrightarrow{\cal A}$ and $\mu':{\cal A}\otimes{\mbox C}\longrightarrow {\cal A}$ are the canonical isomorphisms, defined by
\begin{eqnarray*}
    \mu(k \otimes u) = ku = \mu'(u \otimes k), \qquad \forall u \in {\cal A}, \quad \forall k \in {\mbox C}.
\end{eqnarray*}

\vspace*{0.3cm} {\bf (3)} The C-algebra antihomomorphism (coinverse) $S_{\cal A}: {\cal A} \longrightarrow {\cal A}$ is defined by
\begin{eqnarray*}
    S_{\cal A}(x) & = & x^{-1}, \\ S_{\cal A}(y) & = & - x^{-1} y x^{-1}, \\ S_{\cal A}(z) & = & - z.
\end{eqnarray*}

\vspace{0.3cm}\noindent The coinverse $S$ satisfies
\begin{eqnarray*}
    m \circ (S_{\cal A} \otimes \mbox{id}) \circ \Delta_{\cal A} = \epsilon_{\cal A} = m \circ (\mbox{id} \otimes S_{\cal A}) \circ \Delta_{\cal A}
\end{eqnarray*}
where $m$ stands for the algebra product ${\cal A} \otimes {\cal A} \longrightarrow {\cal A}$.

\vspace{0.3cm}\noindent The coproduct, counit and coinverse which are specified above supply the algebra ${\cal A}$ with a Hopf algebra structure.

\vspace*{0.3cm} \noindent {\bf 2.3 Differential algebra}

\vspace*{0.3cm} \noindent We first note that the properties of the exterior differential {\sf d}. The exterior differential {\sf d} is an operator which gives the mapping from the generators of ${\cal A}$ to the differentials:
\begin{eqnarray*}
    {\sf d} : u \longrightarrow {\sf d}u, \qquad u \in \{x,y,z\}.
\end{eqnarray*}
We demand that the exterior differential {\sf d} has to satisfy two properties: the nilpotency
\begin{eqnarray*}
    {\sf d}^2 & = & 0
\end{eqnarray*}
and the Leibniz rule
\begin{eqnarray*}
    {\sf d}(f g) = ({\sf d} f) g + (-1)^{deg(f)} \, f ({\sf d} g).
\end{eqnarray*}
\noindent A deformed differential calculus on the quantum 3d space is as follows:

\noindent the commutation relations with the coordinates of differentials
\begin{eqnarray*}
  x~ {\sf d} x & = & {\sf d}x~ x, \quad x~ {\sf d}y = q {\sf d}y~ x,  \quad x~ {\sf d}z = q {\sf d}z~ x, \\
  y~ {\sf d}x &=& q^{-1} {\sf d} x~ y, \quad y~ {\sf d} y= {\sf d} y~ y, \quad y~ {\sf d} z = q {\sf d} z~ y,  \\
  z~ {\sf d}x & = & q^{-1}{\sf d}x~ z, \quad z~ {\sf d}y = q^{-1} {\sf d}y~ z, \quad z~ {\sf d}z = {\sf d}z~ z.
\end{eqnarray*}
This algebra is denoted by $\Gamma_1$.

\noindent The commutation relations between the differentials
\begin{eqnarray*}
  {\sf d}x \wedge {\sf d}x & = & 0, \quad {\sf d}y \wedge {\sf d}y  = 0, \quad {\sf d}z \wedge {\sf d}z = 0. \\
  {\sf d}x \wedge {\sf d}y & = & - q~{\sf d}y \wedge {\sf d}x, \\
  {\sf d}y \wedge {\sf d}z & = & - q~{\sf d}z \wedge {\sf d}y, \\
  {\sf d}x \wedge {\sf d}z & = & - q~{\sf d}z \wedge {\sf d}x.
\end{eqnarray*}
This algebra is denoted by $\Gamma_2$.

\vspace*{0.3cm} \noindent
A differential algebra on an associative algebra ${\cal A}$ is a graded associative algebra $\Gamma$ equipped with an operator {\sf d} that has the above properties. Furthermore, the algebra $\Gamma$ has to be generated by $\Gamma_0 \cup \Gamma_1 \cup \Gamma_2$, where $\Gamma_0$ is isomorphic to ${\cal A}$. Let $\Gamma$ be the quoitent algebra of the free associative algebra on the set $\{x,y,z, {\sf d}x, {\sf d}y, {\sf d}z\}$ modulo the ideal $J$ that is generated by the relations of ${\cal R}^3_q$, $\Gamma_1$ and $\Gamma_2$.

\vspace*{0.3cm} \noindent To proceed, one can obtain the relations of the coordinates with their partial derivatives using the expression
\begin{eqnarray*}
    {\sf d} f & = & ({\sf d} x ~\partial_x + {\sf d} y ~\partial_y + {\sf d} z ~\partial_z) f.
\end{eqnarray*}
Consequently one has
\begin{eqnarray*}
\partial_x x &=&1+x \partial_x, \quad \partial_x y = q^{-1} y \partial_x, \quad \partial_x z = q^{-1} z \partial_x, \\
\partial_y x &=& q x \partial_y, \quad \partial_y y = 1+y \partial_y, \quad \partial_y z = q^{-1} z \partial_y, \\
\partial_z x &=& q x \partial_z, \quad \partial_z y = q y \partial_z, \quad \partial_z z = 1 + z \partial_z.
\end{eqnarray*}
Using the fact that ${\sf d}^2 = 0$, one finds
\begin{eqnarray*}
    \partial_x \partial_y & = & q \partial_y \partial_x, \\ \partial_x \partial_z & = & q \partial_z \partial_x, \\
  \partial_y \partial_z & = & q \partial_z \partial_y.
\end{eqnarray*}

\noindent The relations between partial derivatives and differentials are found as
\begin{eqnarray*}
\partial_x {\sf d}x & = & {\sf d} x ~\partial_x, \quad\partial_x {\sf d}y = q^{-1}{\sf d} y ~\partial_x,
   \quad \partial_x {\sf d}z = q^{-1} {\sf d}z \partial_x,  \\
\partial_y {\sf d}x & = & q {\sf d} x ~\partial_y, \quad\partial_y {\sf d}y = {\sf d} y ~\partial_y,
   \quad \partial_y {\sf d}z = q^{-1} {\sf d}z \partial_y, \\
\partial_z {\sf d}x & = & q {\sf d} x ~\partial_z, \quad\partial_z {\sf d}y = q {\sf d} y ~\partial_z,
   \quad \partial_z {\sf d}z = {\sf d}z \partial_z.
\end{eqnarray*}

\noindent We can define three one-forms using the generators of ${\cal A}$. If we call them $\omega_x$, $\omega_y $ and $\omega_z$ then one can define them as follows:
\begin{eqnarray*}
\omega_x & = & {\sf d} x ~x^{-1}, \\
\omega_y & = & {\sf d} y ~x^{-1} - {\sf d} x ~x^{-1}y x^{-1}, \\
\omega_z & = & {\sf d} z.
\end{eqnarray*}
We denote the algebra of forms generated by three elements $\omega_x$, $\omega_y$ and $\omega_z$ by $\Omega$. The
generators of the algebra $\Omega$ with the generators of $\cal A$ satisfy the following rules
\begin{eqnarray*}
x \omega_x & = & \omega_x x, \quad x \omega_y = q \omega_y x, \quad x \omega_z = q \omega_z x, \\
y \omega_x & = & \omega_x y, \quad y \omega_y = q \omega_y y, \quad y \omega_z = q \omega_z y, \\
z \omega_x & = & \omega_x z, \quad z \omega_y = \omega_y z, \quad z \omega_z = \omega_z z.
\end{eqnarray*}
The commutation rules of the generators of $\Omega$ are
\begin{eqnarray*}
\omega_x \wedge \omega_x & = & 0, \quad \omega_y \wedge \omega_y = 0, \quad \omega_z \wedge \omega_z = 0, \\
\omega_x \wedge \omega_y & = & - \omega_y \wedge \omega_x, \\
\omega_y \wedge \omega_z & = & - \omega_z \wedge \omega_y, \\
\omega_x \wedge \omega_z & = & - \omega_z \wedge \omega_x.
\end{eqnarray*}
The algebra $\Omega$ is a graded Hopf algebra [2].

\vspace*{0.3cm} \noindent {\bf 2.4 Lie algebra}

\vspace*{0.3cm} \noindent The commutation relations of Cartan-Maurer forms allow us to construct the algebra of the
generators. In order to obtain the quantum Lie algebra of the algebra generators we first write the Cartan-Maurer forms as
\begin{eqnarray*}
{\sf d} x & = & \omega_x \, x, \\
{\sf d} y & = & \omega_x \, y + \omega_y \, x, \\
{\sf d} z & = & \omega_z.
\end{eqnarray*}
The differantial ${\sf d}$ can then the expressed in the form
\begin{eqnarray*}
    {\sf d} f & = & \left(\omega_x \, T_x + \omega_y \, T_y + \omega_z \, T_z\right) f.
\end{eqnarray*}

\noindent Here $T_x$, $T_y$ and $T_z$ are the quantum Lie algebra generators. Considering an arbitrary function $f$ of the coordinates of the quantum 3d space and using that ${\sf d}^2 = 0$, we find the following commutation relations for the (undeformed) Lie algebra [2]:
\begin{eqnarray*}
    \left[T_x, T_y \right] & = & 0, \\ \left[T_x, T_z \right] & = & 0, \\ \left[T_y, T_z \right] & = & 0.
\end{eqnarray*}

\noindent The commutation relations between the generators of algebra and the coordinates are
\begin{eqnarray*}
T_x \, x & = & x + x \, T_x, \quad T_x \, y = y + y \, T_x, \quad T_x \, z = z\, T_x,  \\
T_y \, x & = & q x \, T_y, \quad T_y \, y = x + q y \, T_y, \quad T_y \,z = z \, T_y, \\
T_z \, x & = & q x \, T_z, \quad T_z \, y = q y \, T_z, \quad T_z \, z = 1 + z \, T_z.
\end{eqnarray*}

\noindent The (quantum) Lie algebra generators can be expressed in terms of the generators of the quantum 3d space and partial differentials:
\begin{eqnarray*}
    T_x & \equiv & x \partial_x + y \partial_y, \\ T_y & \equiv & x \partial_y, \\ T_z & \equiv & \partial_z.
\end{eqnarray*}

\noindent The commutation relations of the Lie algebra generators $T_x$, $T_y$ and $T_z$ with the differentials are following
\begin{eqnarray*}
T_x \, {\sf d} x &=& {\sf d} x \, T_x, \quad T_x \, {\sf d} y={\sf d} y \, T_x, \quad T_x \, {\sf d} z={\sf d}z \,        T_x,\\
T_y \, {\sf d} x & = & q {\sf d}x \, T_y, \quad T_y \, {\sf d} y = q {\sf d} y \, T_y, \quad T_y \, {\sf d} z =
  {\sf d} z \, T_y, \\
T_z \, {\sf d} x & = & q {\sf d} x \, T_z, \quad T_z {\sf d}y=q {\sf d} y T_z, \quad T_z {\sf d}z={\sf d} z \, T_z.
\end{eqnarray*}

\noindent The commutation rules of the Lie algebra generators with one-forms as follows
\begin{eqnarray*}
T_x \, \omega_x &=& \omega_x \, T_x - \omega_x, \quad T_x \, \omega_y = \omega_y \, T_x - \omega_y, \quad T_x \omega_z = \omega_z T_x,  \\
  T_y \, \omega_x & = & \omega_x T_y, \quad T_y \, \omega_y = \omega_y \, T_y -\omega_x, \quad T_y \omega_z =
     \omega_z T_y, \\
  T_z \, \omega_x & = & \omega_x \, T_z, \quad T_z \, \omega_y = \omega_y T_z, \quad T_z \, \omega_z = \omega_z \, T_z.
\end{eqnarray*}

\noindent The Hopf algebra structure of the Lie algebra generators is given by
\begin{eqnarray*}
\Delta(T_x) & = & T_x \otimes {\bf 1} + {\bf 1} \otimes T_x,\\
\Delta(T_y) & = & T_y \otimes {\bf 1} + q^{T_x} \otimes T_y, \\
\Delta(T_z) & = & T_z \otimes {\bf 1} + q^{T_x} \otimes T_z, \\
\epsilon(T_x) & = & 0, \quad \epsilon(T_y) = 0, \quad \epsilon(T_z) = 0, \\
S(T_x) & = & - T_x, \quad S(T_y) = - q^{-T_x} T_y, \quad S(T_z) =
- q^{-T_x} T_z.
\end{eqnarray*}

\vspace*{0.3cm}\noindent
{\bf 2.5 The dual of the Hopf algebra $\cal A$}

\noindent
In this section, in order to obtain the dual of the Hopf algebra ${\cal A}$
defined in section 2, we shall use the method of Refs. 4 and 9.

A pairing between two vector spaces ${\cal U}$ and ${\cal A}$ is a bilinear
mapping
\begin{eqnarray*}
    <,> : {\cal U} \mbox{x} {\cal A} \longrightarrow {\cal C}, \quad
  (u,a) \, \mapsto \,\, <u,a>.
\end{eqnarray*}
We say that the pairing is non-degenerate if
\begin{eqnarray*}
    <u,a> = 0 \quad (\forall a \in {\cal A}) \,\, \Longrightarrow \,\, u = 0
\end{eqnarray*}
and
\begin{eqnarray*}
    <u,a> = 0 \quad (\forall u \in {\cal U}) \,\, \Longrightarrow \,\, a = 0.
\end{eqnarray*}
Such a pairing can be extended to a pairing of ${\cal U} \otimes {\cal U}$
and ${\cal A} \otimes {\cal A}$ by
\begin{eqnarray*}
    <u \otimes v, a \otimes b> = <u,a> <v,b>.
\end{eqnarray*}

Given bialgebras ${\cal U}$ and ${\cal A}$ and a non-degenerate pairing
\begin{eqnarray*}
    <,> : {\cal U} \mbox{x} {\cal A} \longrightarrow {\cal C}  \quad (u,a) \mapsto <u,a> \quad \forall u \in {\cal U}       \quad \forall a \in {\cal A}
\end{eqnarray*}
we say that the bilinear form realizes a duality between ${\cal U}$ and
${\cal A}$, or that the bialgebras ${\cal U}$ and ${\cal A}$ are in duality,
if we have
\begin{eqnarray*}
  <uv, a> & = & <u \otimes v, \Delta_{\cal A}(a)>, \\
  <u, ab> & = & <\Delta_{\cal U}(u), a \otimes b>, \\
  <1_{\cal U}, a> & = & \epsilon_{\cal A}(a),
\end{eqnarray*}
and
\begin{eqnarray*}
    <u, 1_{\cal A}> = \epsilon_{\cal U}(u)
\end{eqnarray*}
for all $u, v \in {\cal U}$ and $a, b \in {\cal A}$.

If, in addition, ${\cal U}$ and ${\cal A}$ are Hopf algebras with coinverse
$\kappa$, then they are said to be in duality if the underlying bialgebras are
in duality and if, moreover, we have
\begin{eqnarray*}
    <S_{\cal U}(u), a> = <u, S_{\cal A}(a)> \quad \forall u \in {\cal U} \quad a \in {\cal A}.
\end{eqnarray*}
It is enough to define the pairing between the generating elements of the
two algebras. Pairing for any other elements of ${\cal U}$ and ${\cal A}$
follows from above relations and the bilinear form inherited by the tensor
product. For example, for
\begin{eqnarray*}
    \Delta_{\cal U}(u) = \sum_k u_k' \otimes u_k'',
\end{eqnarray*}
we have
\begin{eqnarray*}
    <u, ab> = <\Delta_{\cal U}(u), a \otimes b> =
  \sum_k <u_k', a> <u_k'', b>
\end{eqnarray*}
As a Hopf algebra ${\cal A}$ is generated by the elements $x$, $y$ and $z$, and a basis
is given by all monomials of the form
\begin{eqnarray*}
    f = x^k y^l z^m
\end{eqnarray*}
where $k, l, m \in {\cal Z}_+$. Let us denote the dual algebra by
${\cal U}_q$ and its generating elements by $A$ and $B$.

The pairing is defined through the tangent vectors as follows
\begin{eqnarray*}
<X, f> & = & k \, \delta_{l,0} \delta_{m,0}, \\
<Y, f> & = & \delta_{l,1} \delta_{m,0}, \\
<Z, f> & = & \delta_{l,0} \delta_{m,1}.
\end{eqnarray*}
We also have
\begin{eqnarray*}
 <1_{\cal U}, f> = \epsilon_{\cal A}(f) = \delta_{k,0}.
\end{eqnarray*}
Using the defining relations one gets
\begin{eqnarray*}
    <XY, f> = \delta_{l,1} \delta_{m,0}
\end{eqnarray*}
and
\begin{eqnarray*}
    <YX, f> = \delta_{l,1} \delta_{m,0}
\end{eqnarray*}
where differentiation is from the right as this is most suitable for
differentiation in this basis. Thus one obtains one of the commutation relations in
the algebra ${\cal U}_q$ dual to ${\cal A}$ as:
\begin{eqnarray*}
    XY = YX.
\end{eqnarray*}
Similarly, one has
\begin{eqnarray*}
    XZ & = & ZX, \\
    YZ & = & ZY.
\end{eqnarray*}
The Hopf algebra structure of this algebra can be deduced by using the
duality. The coproduct of the elements of the dual algebra is given by
\begin{eqnarray*}
    \Delta_{\cal U}(X) & = & X \otimes 1_{\cal U} + 1_{\cal U} \otimes X, \\
  \Delta_{\cal U}(Y) & = & Y \otimes q^{-X} + 1_{\cal U} \otimes Y, \\
  \Delta_{\cal U}(Z) & = & Z \otimes q^{-X} + 1_{\cal U} \otimes Z.
\end{eqnarray*}
The counity is given by
\begin{eqnarray*}
    \epsilon_{\cal U}(X) = 0, \quad \epsilon_{\cal U}(Y) = 0, \quad \epsilon_{\cal U}(Z) = 0.
\end{eqnarray*}
The coinverse is given as
\begin{eqnarray*}
    S_{\cal U}(X) = - X, \quad S_{\cal U}(Y) = - Y q^X, \quad S_{\cal U}(Z) = - Z q^X.
\end{eqnarray*}

We can now transform this algebra to the form obtained in section 5 by making
the following identities: 
\begin{eqnarray*}
    T_x \equiv X, \quad T_y \equiv q^{X/2} Y q^{X/2}, \quad T_z \equiv q^{X/2} Z q^{X/2}
\end{eqnarray*}
which are consistent with the commutation relation and the Hopf structures.

\newpage\noindent {\bf 3. EXTENDED CALCULUS ON THE QUANTUM 3D SPACE}

\vspace*{0.3cm} \noindent A Lie derivative is a derivation on the algebra of tensor fields over a manifold. The Lie derivative should be defined three ways: on scalar functions, vector fields and tensors.

\noindent The Lie derivative can also be defined on differential forms. In this case, it is closely related to the  exterior derivative. The exterior derivative and the Lie derivative are set to cover the idea of a derivative in  different ways. These differences can be hasped together by introducing the idea of an antiderivation which is called  an inner derivation.

\vspace*{0.3cm} \noindent {\bf 3.1 Inner derivations}

\vspace*{0.3cm} \noindent In order to obtain the commutation rules of the coordinates with inner derivations, we shall use the approach of Ref. 1. Similarly other relations can also obtain. Consequently, we have the following commutation relations:
\begin{itemize}

\item the commutation relations of the inner derivations with $x$, $y$ and $z$
\begin{eqnarray*}
{\bf \textit i}_x \, x &=& x {\bf \textit i}_x, \quad {\bf \textit i}_x \, y=q^{-1} \, y \, {\bf \textit i}_x, \quad       {\bf \textit i}_x \, z=q^{-1} \, z \, {\bf \textit i}_x,\\
{\bf \textit i}_y \, x &=& q \, x \, {\bf \textit i}_y, \quad {\bf \textit i}_y \, y=y \, {\bf \textit i}_y, \quad         {\bf \textit i}_y \, z=q^{-1} \, z \, {\bf \textit i}_y,  \\
{\bf \textit i}_z \, x &=& q \, x \, {\bf \textit i}_z, \quad {\bf \textit i}_z \, y=q \, y \, {\bf \textit i}_z, \quad    {\bf \textit i}_z \, z=z \,{\bf \textit i}_z.
\end{eqnarray*}

\item the relations of the inner derivations with the partial derivatives $\partial_x$, $\partial_y$ and $\partial_z$
\begin{eqnarray*}
{\bf \textit i}_x \, \partial_x &=& \partial_x \, {\bf \textit i}_x, \quad {\bf \textit i}_x \, \partial_y=q \,             \partial_y \, {\bf \textit i}_x, \quad {\bf \textit i}_x \, \partial_z=q \, \partial_z \, {\bf \textit i}_x,\\
{\bf \textit i}_y \, \partial_x &=& q^{-1} \, \partial_x \, {\bf \textit i}_y, \quad {\bf \textit i}_y \, \partial_y=       \, \partial_y \, {\bf \textit i}_y, \quad {\bf \textit i}_y \, \partial_z=q \, \partial_z \, {\bf \textit  i}_y,  \\
{\bf \textit i}_z \, \partial_x &=& q^{-1} \, \partial_x \, {\bf \textit i}_z , \quad {\bf \textit i}_z \, \partial_y      =q^{-1} \, \partial_y \, {\bf \textit i}_z, \quad {\bf \textit i}_z \, \partial_z=\partial_z \, {\bf \textit i}_z.
\end{eqnarray*}

\item the commutation relations between the differentials and the inner derivations
\begin{eqnarray*}
{\bf \textit i}_x \wedge \, {\sf d}x &=& 1-{\sf d}x \wedge \, {\bf \textit i}_x, \quad {\bf \textit i}_x \wedge \, {\sf     d}y=-q^{-1} \, {\sf d}y \wedge \, {\bf \textit i}_x, \\
{\bf \textit i}_y \wedge \, {\sf d}x &=& -q \, {\sf d}x \wedge \, {\bf \textit i}_y, \quad {\bf \textit i}_y \wedge \,     {\sf d}y=1-{\sf d}y \wedge \, {\bf \textit i}_y, \\
{\bf \textit i}_z \wedge \, {\sf d}x &=& -q \, {\sf d} x \wedge \, {\bf \textit i}_z, \quad {\bf \textit i}_z \wedge \,    {\sf d}y=-q \, {\sf d}y \wedge \, {\bf \textit i}_z, \\
{\bf \textit i}_x \wedge \, {\sf d}z & = & -q^{-1} \, {\sf d}      z \wedge \, {\bf \textit i}_x, \quad {\bf \textit i}_y \wedge \, {\sf d}z=-q^{-1} \, {\sf d}z \wedge \, {\bf \textit i}_y, \\
{\bf \textit i}_z\wedge \, {\sf d}z & = & 1-{\sf d} z \wedge \,    {\bf \textit i}_z.
\end{eqnarray*}
\end{itemize}

\noindent {\bf 3.2 Lie derivations}

\vspace*{0.3cm} \noindent In this section we find the commutation rules of the Lie derivatives with functions, {\it i.e.} the elements of the algebra ${\cal A}$, their differentials, etc., using the approach of [1] as follows:
\begin{itemize}
    \item the relations between the Lie derivatives and the elements of ${\cal A}$
\begin{eqnarray*}
{\cal L}_x\, x &=&1+x\, {\cal L}_x, \quad {\cal L}_x \, y=q^{-1} \, y \, {\cal L}_x, \quad {\cal L}_x \, z=q^{-1} \, z    \, {\cal L}_x ,\\
{\cal L}_y \, x &=&q \, x \, {\cal L}_y, \quad {\cal L}_y \, y=1+y \, {\cal L}_y, \quad {\cal L}_y \, z=q^{-1} \, z \,    {\cal L}_y, \\
{\cal L}_z \, x &=&q \, x \, {\cal L}_z, \quad {\cal L}_z \, y=q \, y \, {\cal L}_z, \quad {\cal L}_z \, z=1+z \, {\cal   L}_z .
\end{eqnarray*}
\item The relations of the Lie derivatives with the differentials
\begin{eqnarray*}
{\cal L}_x \, {\sf d} x &=& {\sf d}x \, {\cal L}_x, \quad {\cal L}_x \, {\sf d}y=q^{-1} \, {\sf d}y \, {\cal L}_x,          \quad {\cal L}_x \, {\sf d}z=q^{-1} \, {\sf d}z \, {\cal L}_x,\\
{\cal L}_y \, {\sf d} x &=&q \, {\sf d}x \, {\cal L}_y, \quad {\cal L}_y \, {\sf d}y={\sf d}y \, {\cal L}_y, \quad          {\cal L}_y \, {\sf d}z=q^{-1} \, {\sf d}z \, {\cal L}_y, \\
{\cal L}_z \, {\sf d} x &=&q \, {\sf d}x \, {\cal L}_z, \quad {\cal L}_z \, {\sf d}y=q \, {\sf d}y \, {\cal L}_z, \quad     {\cal L}_z \, {\sf d}z={\sf d}z \, {\cal L}_z .
\end{eqnarray*}
\end{itemize}

\noindent Other commutation relations can be similarly obtained. To complete the description of the above scheme, we get below the remaining commutation relations as follows:
\begin{itemize}
    \item the Lie derivatives and partial derivatives
\begin{eqnarray*}
   {\cal L}_x\, \partial_x & = & \, \partial_x\, {\cal L}_x ,\quad \ {\cal L}_x \, \partial_y  =  q \, \partial_y\, {\cal L}_x ,\quad \ {\cal L}_x \, \partial_z   =  q \,\partial_z  \, {\cal L}_x,\\
   {\cal L}_y\, \partial_x & = & \, \, q^{-1} \,\partial_x  \, {\cal L}_y ,\quad \ {\cal L}_y \, \partial_y  = \partial_y\, {\cal L}_y ,\quad \ {\cal L}_y \, \partial_z   =  q \,\partial_z  \, {\cal L}_y , \\
   {\cal L}_z\, \partial_x & = & \, \, q^{-1} \,\partial_x  \, {\cal L}_z ,\quad \ {\cal L}_z \, \partial_y  =  q^{-1} \, \partial_y\, {\cal L}_z ,\quad \ {\cal L}_z \, \partial_z   =  \partial_z  \, {\cal L}_z .
\end{eqnarray*}
  \item the inner derivations
\begin{eqnarray*}
{\bf \textit i}_x \wedge \, {\bf \textit i}_y & = & - q \, {\bf \textit i}_y \wedge \, {\bf \textit i}_x,\\
{\bf \textit i}_x \wedge \, {\bf \textit i}_z & = & - q \, {\bf \textit i}_z \wedge \, {\bf \textit i}_x ,  \\
{\bf \textit i}_y \wedge \, {\bf \textit i}_z & = & - q \, {\bf \textit i}_z \wedge \, {\bf \textit i}_y .
\end{eqnarray*}
  \item the Lie derivatives and the inner derivations
\begin{eqnarray*}
{\cal L}_x\, {\bf \textit i}_x &=& {\bf \textit i}_x \, {\cal L}_x, \quad {\cal L}_x \, {\bf \textit i}_y=q \, {\bf         \textit i}_y \, {\cal L}_x, \quad {\cal L}_x \, {\bf \textit i}_z=q \, {\bf \textit i}_z \, {\cal L}_x,\\
{\cal L}_y \, {\bf \textit i}_x &=& q^{-1} \, {\bf \textit i}_x \, {\cal L}_y, \quad {\cal L}_y \, {\bf \textit i}_y=       {\bf \textit i}_y \, {\cal L}_y, \quad {\cal L}_y \, {\bf \textit i}_z=q \, {\bf \textit i}_z \, {\cal L}_y, \\
{\cal L}_z \, {\bf \textit i}_x &=& q^{-1} \, {\bf \textit i}_x  \, {\cal L}_z, \quad {\cal L}_z \, {\bf \textit i}_y=      q^{-1} \, {\bf \textit i}_y \, {\cal L}_z, \quad {\cal L}_z \, {\bf \textit i}_z={\bf \textit i}_z \, {\cal L}_z.
\end{eqnarray*}
  \item the Lie derivatives
\begin{eqnarray*}
{\cal L}_x {\cal L}_y & = & q\, {\cal L}_y {\cal L}_x, \\
{\cal L}_x {\cal L}_z & = & q\, {\cal L}_z {\cal L}_x, \\
{\cal L}_y {\cal L}_z & = & q\, {\cal L}_z {\cal L}_y.
\end{eqnarray*}
\end{itemize}

Note that the Lie derivatives can be written as follows:
\begin{eqnarray*}
{\cal L}_x & = & x^{-1} \, T_x - x^{-1} y x^{-1} \, T_y,\\
{\cal L}_y & = & x^{-1} \, T_y, \\
{\cal L}_z & = & T_z.
\end{eqnarray*}

\noindent {\bf ACKNOWLEDGMENT}

\vspace*{0.3cm}\noindent This work was supported in part by {\bf TBTAK} the Turkish Scientific and Technical Research Council.

\newpage\noindent {\bf REFERENCES}
\begin{itemize}
\item [1.] Celik, Salih: {\it J. Math. Phys.} {\bf 47 (8)}: Art. No: 083501
\item [2.] Celik, Sultan A. and Yasar, E.: {\it Czech. J. Phys.} {\bf 56} (2006), 229.
\item [3.] Chryssomalakos, C., Schupp P. and Zumino, B.: "Induced extended calculus on the quantum plane", hep-th /9401141.
\item [4.] Dobrev, V. K.: {\it J. Math. Phys.} {\bf 33} (1992), 3419.
\item [5.] Manin, Yu I.: "Quantum groups and noncommutative geometry", (Montreal Univ. Preprint, 1988).
\item [6.] Schupp, P., Watts, P., Zumino, B.: {\it Lett. Math. Phys.} {\bf 25} (1992), 139.
\item [7.] Schupp, P., Watts P., Zumino, B.: "Cartan calculus on quantum Lie algebras", hep-th/9312073.
\item [8.] Schupp, P.: "Cartan calculus: Differential geometry for quantum groups", hep-th/9408170.
\item [9.] A. Sudbery, A.: {\it Proc. Workshop on Quantum Groups}, Argogne (1990) eds. T. Curtright, D. Fairlie and C. Zachos, pp. 33-51.
\item [10.] Wess, J. and Zumino, B.: {\it Nucl. Phys.} (Proc. Suppl.) {\bf 18B} (1990), 302.
\item [11.] Woronowicz, S.L.: {\it Commun. Math. Phys.} {\bf 111} (1987), 613.
\item [12.] Woronowicz, S.L.  {\it Commun. Math. Phys.} {\bf 122} (1989), 125.
\end{itemize}

\end{document}